\documentclass[12pt, a4paper,reqno]{amsart}
\usepackage{amscd,amssymb,amsmath, array,latexsym,amsbsy}
\input{diagrams.tex}
\allowdisplaybreaks[2]
\begin{document}

\def\Ind{\operatorname{Ind}}
\def\Hom{\operatorname{Hom}}
\def\Aut{\operatorname{Aut}}
\def\Rep{\operatorname{Rep}}
\def\Ad{\operatorname{Ad}}
\def\sp{\operatorname{sp}}
\def\dashind{\operatorname{\!-Ind}}
\def\id{\operatorname{id}}
\def\rt{\operatorname{{rt}}}
\def\lt{\operatorname{{lt}}}
\def\res{\operatorname{Res}}

\def\H{\mathcal{H}}
\def\A{\mathcal{A}}
\def\C{\mathbb{C}}
\def\T{\mathbb{T}}
\def\Z{\mathbb{Z}}
\def\R{\mathbb{R}}

\newtheorem{thm}{Theorem}
\newtheorem{cor}[thm]{Corollary}
\theoremstyle{definition}
\newtheorem{defn}[thm]{Definition}
\newtheorem{remark}[thm]{Remark}

\title[\boldmath Extending unitary representations of  
subgroups]{Twisted actions and the obstruction to extending\\  unitary   
representations of subgroups}

\author[an Huef]{Astrid an Huef}
\address{School of Mathematics\\
University of New South Wales\\
Sydney, NSW 2052\\
Australia}
\email{astrid@unsw.edu.au}

\author[Raeburn]{Iain Raeburn}
\address{School of Mathematical and Physical Sciences\\
University of Newcastle\\
NSW 2308\\ Australia}
\email{iain@frey.newcastle.edu.au}


\date{September 22, 2003}

\begin{abstract} Suppose that $G$ is a locally compact group and  $\pi$  
is a  (not necessarily irreducible) unitary representation of a closed  
normal subgroup $N$ of $G$ on a Hilbert space $\H$. We extend results  
of Clifford and Mackey to determine when $\pi$ extends to a unitary  
representation  of $G$ on the same space $\H$ in terms of a  
cohomological obstruction.
\end{abstract}

\thanks{This research was supported by grants from the Australian  
Research Council and the University of New South Wales.}
\maketitle

Let $G$ be a  group and $\pi:N\to U(\H)$ a unitary representation of a  
normal
subgroup $N$ of $G$.  When is $\pi$ the restriction of a unitary  
representation of $G$?

If $\pi$ does extend to a representation $\rho$ of $G$, then $\rho(s)$  
implements a unitary equivalence between $\pi$ and  
$\pi^s:n\mapsto\pi(sns^{-1})$.  So an obvious necessary condition is  
that $\pi$ should be equivalent to $\pi^s$ for each $s\in G$ (we say  
that $\pi$ is \emph{$G$-invariant}), and the problem is to decide when  
a $G$-invariant representation extends.

Clifford answered this extension problem in \cite{Clifford} when $G$ is  
discrete,  $\pi$ is irreducible and $\H$ is finite-dimensional. In  
modern language, Clifford showed that if $\pi$ is $G$-invariant, then  
there is an obstruction to extending the representation in the  
cohomology group  $H^2(G/N,\T)$, where $\T$ is the unit  circle.
Mackey  extended Clifford's result to irreducible unitary  
representations of normal closed subgroups of locally compact groups  
\cite[Theorem~8.2]{Mackey}.  Mackey's solution involves Borel cocycles,  
so his obstruction lies in a cohomology theory where  all the cochains  
are Borel.
The resulting cohomology groups were subsequently analysed by Moore in  
\cite{Moore1,Moore2,Moore3,Moore4}.

This extension problem has recently resurfaced in the context of  
compact Lie groups in \cite{CKS}, where it was tackled using  the  
structure theory of Lie groups, and in \cite{BE}, where it was studied  
in the context of nonabelian duality for locally compact groups and  
crossed products of $C^*$-algebras.
Here we investigate a cohomological obstruction to the extension of an  
arbitrary $G$-invariant unitary representation $\pi$ of $N$, and its  
relationship to the results in \cite{CKS} and \cite{BE}.  Our  
obstruction is a twisted action of $G/N$ on the von Neumann algebra  
$\pi(N)'$ of operators which commute with every $\pi(n)$; the  
representation extends if and only if this twisted action is  
equivalent, in a natural sense, to an ordinary action. We then use a  
stabilisation trick to show that if $\pi$ is $G$-invariant then  
infinite multiples $\pi\otimes 1$ of $\pi$ always extend.

\subsection*{Preliminaries}
Let $G$ be a second-countable locally compact group with a closed  
normal subgroup $N$.
  We endow the group $U(\H)$ of all unitary operators on a separable  
Hilbert space $\H$  with the strong operator topology, and note that  
$U(\H)$ is a Polish group (in the sense that the topology is given by a  
complete metric). A unitary representation $\rho$ of $G$ is a  
continuous  homomorphism $\rho:G\to U(\H)$.  A function $f:G\to\H$ is  
Borel if $f^{-1}(O)$ is a Borel set for each open set $O$ of $\H$;  
equivalently, if $s\mapsto (f(s)\, |\, h):G\to \C$ is a Borel function  
for each $h\in\H$. We use a left-invariant Haar measure on $G$.

Let $\A$ be a von Neumann algebra acting on a separable Hilbert space  
$\H$. The group $U(\A)$ of unitary elements in $\A$  is a Polish group 
in the ultra-weak topology,
and it is then a closed subgroup of $U(\H)$.  
The group $\Aut(\A)$ 
of automorphisms of $\A$  is Polish in the topology of pointwise ultra-weak convergence;  this is called the  
$u$-topology in \cite[Definition~3.4]{Haagerup}.  For $u\in U(\A)$, we denote by $\Ad u$  
the automorphism $a\mapsto uau^*$ of $\A$, and note that $\Ad:U(\A)\to  
\Aut(\A)$ is a continuous homomorphism.

\begin{defn}
A \emph{twisted action} of $G$ on  $\A$ is a pair $(\alpha, \sigma)$ of  
maps $\alpha:G\to \Aut (\A)$ and $\sigma:G\times G\to U(\A)$ such that
\begin{enumerate}
\item $\alpha$ and $\sigma$ are Borel,
\item $\alpha_e=\id, \sigma(e,s)=\sigma(s,e)=1$ for $s\in G$,
\item $\alpha_s\circ\alpha_t=\Ad\sigma(s,t)\circ\alpha_{st}$ for  
$s,t\in G$, and
\item $\alpha_r(\sigma(s,t))\sigma(r,st)=\sigma(r,s)\sigma(rs,t)$ for  
$r,s,t\in G$.
\end{enumerate}
Two twisted actions $(\alpha,\sigma)$ and $(\beta,\omega)$ of $G$ on  
$\A$ are \emph{exterior equivalent} if there is a Borel map $\nu:G\to  
U(\A)$ such that
\begin{enumerate}
\item $\beta_s=\Ad \nu_s\circ\alpha_s$, and
\item $\omega(s,t)=\nu_s\alpha_s(\nu_t)\sigma(s,t)\nu_{st}^*$.
\end{enumerate}
\end{defn}

These definitions are the von-Neumann algebraic analogues of  
\cite[Definitions~2.1 and~3.1]{PR}.
Our definition of twisted action is slightly  different from the one  
used in \cite[Definition~2.1]{Colin1}, where the map $s\mapsto\alpha_s$  
is required to be continuous.

\subsection*{Main Results} In Theorem~\ref{thm-main} we prove that the  
obstruction to extending a $G$-invariant unitary representation $\pi$  
of $N$ is a twisted action of $G/N$ on the von Neumann algebra  
$\pi(N)'$, and in Theorem~\ref{connect}  we discuss the extension  
problem in the context of non-abelian duality for amenable groups $G$.  
We reconcile  the two approaches in Remark~\ref{remark}; to do so one  
needs to understand not only the statement of Theorem~\ref{thm-main}  
but also its proof.

\begin{thm}\label{thm-main}
Let $N$ be a closed normal subgroup of a second-countable locally  
compact group $G$.  Suppose  $\pi:N\to U(\H)$ is a unitary  
representation of $N$ which is $G$-invariant.  Then there is a twisted  
action $(\alpha,\sigma)$  of $G/N$ on the commutant $\pi(N)'$ of  
$\pi(N)$ such that $\pi$ extends to a  unitary representation $\rho$ of  
$G$ on $\H$  if and only if $(\alpha,\sigma)$ is exterior equivalent to  
an action.
\end{thm}

\begin{proof}
We start by constructing the twisted action $(\alpha,\sigma)$.
Since $\pi^s$ is unitarily equivalent to $\pi$ for all $s\in G$,  there  
exist unitary operators $W_s\in U(\H)$ such that  
$W_s\pi(n)W_s^*=\pi(sns^{-1})$.  We claim that we can choose  $W_s$  
such that the map $s\mapsto W_s$ is  Borel.  To see this, let
\[
H=\{(W,s):W\in U(\H), s\in G\text{ \ and\ }  
W\pi(n)W^*=\pi(sns^{-1})\text{\ for $n\in N$}\}.\]
Then $H$ is a subgroup of $U(\H)\times G$; we claim that $H$ is closed.  
So suppose that the net $(W_\beta,s_\beta)\in H$ converges   
to $(W,s)$. Then $W_\beta\pi(n)$ converges strongly to $W\pi(n)$ for  
each $n\in N$, and since   multiplication in $U(\H)$ is jointly  
continuous, $\pi(s_\beta ns_\beta^{-1})W_\beta$ converges strongly to  
$\pi(s ns^{-1})W$.  Thus $W\pi(n)=\pi(s ns^{-1})W$, and $H$ is closed.  
Now $H$ is Polish since it is a  closed subgroup of  a Polish group, and 
the quotient map $H\to G:(W,s)\mapsto s$ has a Borel section $s\mapsto  
(W_s,s)$ by \cite[Proposition~4]{Moore3}.

Another application of \cite[Proposition~4]{Moore3} shows that the  
quotient map $G\to G/N$ admits a Borel section $c:G/N\to G$. We set
\begin{equation}\label{eq-V}
V_s=W_{c(sN)}\pi(c(sN)^{-1}s).
\end{equation}
Then $s\mapsto V_s$ is Borel and $V_{sn}=V_s\pi(n)$ for $s\in G$ and  
$n\in N$.
We define $\alpha_s=\Ad V_s$. For $T\in \pi(N)'$ and $s\in G$, we have
\begin{align*}
\alpha_s(T)\pi(n)&=V_sTV_s^*\pi(n)=V_sT(V_s^*\pi(n)V_s)V_s^*
=V_sT\pi(s^{-1}ns)V_s^*\\&=V_s\pi(s^{-1}ns)TV_s^*
=V_s\pi(s^{-1}ns)V_s^*V_sTV_s^*=\pi(n)\alpha_s(T);
\end{align*}
thus $\alpha_s:\pi(N)'\to\pi(N)'$, and $\alpha_s$ is an automorphism of  
$\pi(N)'$ because $V_s$ is unitary.

To see that $s\mapsto \alpha_s:G\to \Aut(\pi(N)')$ is Borel, we will  
show that if $V_\beta$ converges to $V$ in the strong operator topology  
and $\Ad V_\beta$ and $\Ad V$ leave $\pi(N)'$ invariant, then $\Ad  
V_\beta$ converges to $\Ad V$ in $\Aut(\pi(N)')$. It then follows that  
$s\mapsto \alpha_s:G\to \Aut(\pi(N)')$  is Borel because $s\mapsto V_s$  
is Borel.
The topology on $\Aut(\pi(N)')$ is the topology generated by the  
seminorms $\alpha\mapsto \|f\circ \alpha\|$, where $f\in \pi(N)'_*$   
and the pre-dual $\pi(N)'_*$ has been identified with the ultra-weakly  
continuous functionals on $\pi(N)'$.  The ultra-weakly continuous  
functionals on $\pi(N)'$ have  the form $f(T)=\sum_{n=1}^\infty  
(Th_n\,|\, k_n)$ where $h_n,k_n\in\H$ satisfy  
$\sum_{n=1}^\infty\|h_n\|^2,\sum_{n=1}^\infty\|k_n\|^2 <\infty$ (see,  
for example, \cite[pages 482--483]{KR}). Let $\epsilon>0$.  If $K$  is  
the maximum of $\big(\sum_{n=1}^\infty\|h_n\|^2\big)^{1/2}$ and
$\big(\sum_{n=1}^\infty\|k_n\|^2\big)^{1/2}$, then
\begin{align}
\|f\circ \Ad V_\beta&-f\circ \Ad V\|
=\sup\big\{ \|f(V_\beta T V_\beta^*-VTV^*)\|: \|T\|=1,  
T\in\pi(N)'\big\}\notag
\\
&=\sup\Big\{ \Big| \sum_{n=1}^\infty ( V_\beta^*h_n\, |\,  
T^*V_\beta^*k_n)-(TV^*h_n\, |\, V^*k_n) \Big|:\|T\|=1,  
T\in\pi(N)'\Big\}\notag
\\
&\leq \sup\Big\{  \sum_{n=1}^\infty
\|(V_\beta^*-V^*)h_n\|\,\|T^*V_\beta^*k_n\|+\|TV^*h_n\|\,\|(V_\beta^*- 
V^*)k_n\|\Big\}
\notag\\
&\leq  
\sum_{n=1}^\infty\|(V_\beta^*-V^*)h_n\|\,\|k_n\|+\|h_n\|\,\|(V_\beta^*- 
V^*)k_n\|\notag\\
&\leq K\Big(\sum_{n=1}^\infty\|(V_\beta^*-V^*)h_n\|^2\Big)^{1/2}+
K\Big(\sum_{n=1}^\infty\|(V_\beta^*-V^*)k_n\|^2\Big)^{1/2}
\label{eq-borel}
\end{align}
by H\"older's inequality.
Since each $V_\beta$ is a normal operator, we have  $V_\beta^*\to V^*$. Now choose  
$N>0$ such that   
$\sum_{n=N}^\infty\|h_n\|^2<{\epsilon^2}({16K^2})^{-1}$ and  
$\sum_{n=N}^\infty\|k_n\|^2<{\epsilon^2}({16K^2})^{-1}$. Then,
for each $n<N$, choose an open neighbourhood $O_n$ of $V^*$ such that
\[
\|(V_\beta^*-V^*)h_n\|^2<\frac{\epsilon^2}{8(N- 
1)K^2}\quad\text{and}\quad
\|(V_\beta^*-V^*)k_n\|^2<\frac{\epsilon^2}{8(N-1)K^2}
\]
whenever $V_\beta^*\in O_n$, and check that if   
$V_\beta^*\in\bigcap_{n=1}^{N-1}O_n$ then $\eqref{eq-borel}<\epsilon$.
This proves that $\Ad V_\beta$ converges to $\Ad V$, and it follows that $\alpha:G\to\Aut(\pi(N)'):s\mapsto\Ad  
V_s$ is Borel.

Next we define $\sigma(s,t)=V_sV_tV_{st}^*$.  Then
\begin{align*}
\sigma(s,t)\pi(n)&=V_sV_tV_{st}^*\pi(n)=V_sV_t\pi((st)^{- 
1}nst)V_{st}^*\\
&=V_sV_t\pi((st)^{-1}nst)V_t^*V_tV_{st}^*\\
&=V_s\pi(s^{-1}ns)V_tV_{st}^*
=\pi(n)V_sV_tV_{st}^*
\end{align*}
for all $n\in N$, so $\sigma:G\times G\to U(\pi(N)')$. Note that  
$\sigma$ is Borel because $s\mapsto V_s$ is Borel and both $V_s\mapsto  
V_s^*$ and $(s,t)\mapsto st$ are continuous. The equation  
$V_{sn}=V_s\pi(n)$ implies  that $\sigma(s,n)=1=\sigma(n,s)$ for  $s\in  
G$ and $n\in N$. We have
\begin{align*}
\sigma(r,s)\sigma(rs,t)&=V_rV_sV_{rs}^*V_{rs}V_tV_{rst}^*=V_rV_sV_tV_{rs 
t}^*\\
&=V_rV_sV_t(V_{st}^*V_r^*V_rV_{st})V_{rst}^*=\alpha_r(\sigma(s,t))\sigma 
(r,st),
\end{align*}
and, for $T\in\pi(N)'$,
\begin{equation}\label{eq-twist-alpha}
\alpha_s(\alpha_t(T))=V_sV_tTV_t^*V_s^*=V_sV_tV_{st}^*V_{st}TV_{st}^*V_{ 
st}V_t^*V_s^*
=\sigma(s,t)\alpha_{st}(T)\sigma(s,t)^*.
\end{equation}
Thus $(\alpha,\sigma)$ is a twisted action of $G$ on $\pi(N)'$. But  
$\alpha_s$ depends only on $sN$ since
$$\alpha_{sn}(T)=V_s\pi(n)T\pi(n)^*V_s^*
=V_s\pi(n)\pi(n)^*TV_s^*=\alpha_s(T)$$ for all $n\in N$. We also have
\[
\sigma(s,tn)=V_sV_{tn}V_{stn}^*=V_sV_t\pi(n)^*V_{st}^*=\sigma(s,t),
\]
and $\sigma(sn,t)=\sigma(sns^{-1}s,t)=\sigma(s,t)$ because  
$\sigma(ns,t)=\sigma(s,t)$.
So we can view $(\alpha, \sigma)$ as a twisted action of $G/N$ on  
$\pi(N)'$.

Now
suppose that $\pi$ extends to a continuous representation $\rho$ of $G$  
on $\H$.  Let $V_s$ be as  in \eqref{eq-V} and define $\nu:G/N\to  
U(\H)$ by $\nu_{sN}=\rho(s)V_s^*$.  Then $\nu$ is Borel because  
$s\mapsto V_s^*$ is Borel and $\rho$ is continuous, and
\begin{equation}\label{eq-trivial}
\omega(sN,tN):=\nu_{sN}\alpha_s(\nu_{tN})\sigma(sN,tN)\nu_{stN}^*=1.
\end{equation}
If $\beta_{sN}:=\Ad\nu_{sN}\circ\alpha_{sN}$, then we can deduce from  
\eqref{eq-trivial} and \eqref{eq-twist-alpha} that
\begin{align*}
\beta_{sN}\circ\beta_{tN}
&=(\Ad\nu_{sN}\circ\alpha_{sN})\circ(\Ad\nu_{tN}\circ\alpha_{tN})
=\Ad\nu_{sN}\circ\alpha_{sN}\circ\Ad\nu_{tN}\circ\alpha_{tN}\\
&=\Ad(\nu_{sN}\alpha_{sN}(\nu_{tN}))\circ\alpha_{sN}\circ\alpha_{tN}
=\Ad(\nu_{stN}\sigma(sN,tN)^*)\circ\alpha_{sN}\circ\alpha_{tN}\\
&=\Ad\nu_{stN}\circ\alpha_{stN}=\beta_{stN}.
\end{align*}
Since $\beta:G/N\to \Aut(\pi(N)')$ is a Borel homomorphism between  
Polish groups, it  is continuous by \cite[Proposition~5]{Moore3}. Thus  
$(\beta,1)$ is an ordinary action, and $(\alpha,\sigma)$ is exterior  
equivalent to an action.

Conversely, if $(\alpha,\sigma)$ is exterior equivalent to an action,  
then there exists a Borel map $\nu:G/N\to U(\H)$ such that  
$\nu_{sN}\alpha_{sN}(\nu_{tN})\sigma(sN,tN)\nu_{stN}^*=1$.
Set $\rho(s)=\nu_{sN}V_s$.  Then
\begin{align*}
\rho(s)\rho(t)
&=\nu_{sN}V_s\nu_{tN}V_t=\nu_{sN}\alpha_{sN}(\nu_{tN})V_sV_t\\
&=
\nu_{sN}\alpha_{sN}(\nu_{tN})\sigma(sN,tN)V_{st}
=\nu_{stN}V_{st}=\rho(st).
\end{align*}
Thus $\rho:G\to U(\H)$ is a  Borel homomorphism between Polish groups,  
and hence is continuous by \cite[Proposition~5]{Moore3}; $\rho$ is the  
required extension of $\pi$.
\end{proof}

\begin{cor}\label{cor-main}
If $\pi:N\to U(\H)$ is a $G$-invariant unitary representation of $N$,  
then there is a unitary representation $\rho$ of $G$ on $\H\otimes  
L^2(G/N)$ such that $\rho|_N=\pi\otimes 1$.
\end{cor}

\begin{proof}
Let $(\alpha,\sigma)$ be the twisted action of $G/N$ on $\pi(N)'$  
constructed above.  Then the twisted action  
$(\beta,\omega):=(\alpha\otimes \id,\sigma\otimes 1)$  of $G/N$ on  
$\pi(N)'\otimes B(L^2(G/N))=(\pi\otimes 1)(N)'$ is the obstruction to  
extending $\pi\otimes 1$. We will show that $(\beta,\omega)$ is  
exterior equivalent to an action. Similar ``stabilisation tricks'' have  
been used in \cite[Proposition~2.1.3]{Colin2} and  
\cite[Theorem~3.4]{PR}, for example.

We begin by  identifying $\H\otimes L^2(G/N)$ with the space  
$L^2(G/N,\H)$ of Bochner square-integrable functions. Since $\H$ is  
separable, $\xi\in L^2(G/N,\H)$ if and only if $\xi$ is a Borel  
function from $G/N$ to $\H$ and $\int_{G/N} \|\xi(sN)\|^2\,  
d(sN)<\infty$. Define $\nu:G/N\to U(L^2(G/N,\H))$ by
\[
(\nu_{tN}\xi)(rN)=\sigma(tN,t^{-1}r^{-1}N)^*\xi(rtN)\Delta(tN)^{1/2},
\]
where $\Delta$ is the modular function of $G/N$ and $\xi\in  
L^2(G/N,\H)$. (The modular function is necessary to ensure that  
$\nu_{tN}$ is unitary.)
Then
\[
(\nu_{tN}^*\xi)(rN)=\sigma(tN,r^{-1}N)\xi(rt^{-1}N)\Delta(tN)^{-1/2},
\]
and hence
\begin{align*}
(\beta_{sN}(\nu_{tN}^*)\nu_{sN}^*\nu_{stN}\xi)(rN)
&=
\alpha_{sN}(\sigma(tN,r^{-1}N))(\nu_{sN}^*\nu_{stN}\xi)(rt^{- 
1}N)\Delta(tN)^{-1/2}\\
&=
\alpha_{sN}(\sigma(tN,r^{-1}N))\sigma(sN,tr^{-1}N)
(\nu_{stN}\xi)(rt^{-1}s^{-1}N)\Delta(stN)^{-1/2}\\
&=
\alpha_{sN}(\sigma(tN,r^{-1}N))\sigma(sN,tr^{-1}N)\sigma(stN,r^{- 
1}N)^*\xi(rN)\\
&=\sigma(sN,tN)\xi(rN)\\
&=(\sigma(sN,tN)\otimes 1)\xi(rN)\\
&=\omega(sN,tN)\xi(rN).
\end{align*}
It follows that
\begin{equation}\label{eqtrivial}
\nu_{sN}\beta_{sN}(\nu_{tN})\omega(sN,tN)\nu_{stN}^*=1.
\end{equation}
If we now define $\gamma:G/N\to \Aut(\pi(N)')$ by $\gamma_{sN}=\Ad  
\nu_{sN}\circ\beta_{sN}$, then (\ref{eqtrivial}) implies that $\gamma$  
is a homomorphism. It remains to show that $\nu$ is Borel, and it then  
follows from \cite[Proposition~5]{Moore3} that  
$\gamma=\Ad\nu\circ\beta:G/N\to\Aut(\pi(N)')$ is continuous.

Since $U(L^2(G/N,\H))$ has the strong operator topology,
$\nu$ is Borel if and only if
$sN\mapsto \nu_{sN}\xi$ is Borel for every $\xi\in L^2(G/N,\H)$, and  
hence
if and only if $sN\mapsto(\nu_{sN}\xi\,|\,\eta)$ is Borel for every  
$\xi,\eta\in L^2(G/N,\H)$.  Since $(U,h)\mapsto Uh$ is continuous, the  
map $(sN,tN,rN)\mapsto (\sigma(sN,tN),\xi(rN))\mapsto  
\sigma(sN,tN)\xi(rN)$ is Borel, and hence so is
\begin{equation}\label{eq1}
(tN,rN)\mapsto  
\big|\big(\sigma(tN,t^{-1}r^{-1}N)^*\xi(rtN)\,\big|\,\eta(rN)\big)\big|
\end{equation}
Since \eqref{eq1} is dominated by $\|\xi(rtN)\|\,\|\eta(rN)\|$, and an  
application of Tonelli's Theorem shows that this is integrable over  
$G/N\times G/N$,
it follows from  Fubini's Theorem  that
\begin{equation*}
tN\mapsto  
\int_{G/N}\big(\sigma(tN,t^{-1}r^{- 
1}N)^*\xi(rtN)\,\big|\,\eta(rN)\big)\, d(rN)
\end{equation*}
defines, almost everywhere, an integrable (and therefore Borel)  
function.
Multiplying by $\Delta(tN)^{1/2}$ shows that   
$tN\mapsto(\nu_{tN}\xi\,|\,\eta)$ is Borel. Thus $\nu$ is Borel and  
$\gamma$ is continuous.

Thus $\nu$ implements an exterior equivalence between $(\beta,\omega)$  
and the ordinary action $(\gamma, 1)$. It now follows from  
Theorem~\ref{thm-main} that there is a representation $\rho$ of $G$  
with $\rho|_N=\pi\otimes 1$.
\end{proof}

\subsection*{The irreducible case} When the representation $\pi$ of $N$  
is irreducible, the commutant  $\pi(N)'$ is $\C1$, the action $\alpha$  
is trivial, and the obstruction $\sigma$ to extending $\pi$ is a Borel  
cocycle in the Moore cohomology group $H^2(G/N,\T)$. Thus we recover  
Mackey's \cite[Theorem~8.2]{Mackey} as it applies to ordinary (that is,  
non-projective) irreducible representations.

When the obstruction $\sigma$ is non-trivial, we can recover  
Corollary~\ref{cor-main} from another important part of the Mackey  
machine \cite[Theorem~8.3]{Mackey}:
$\pi$ extends to a projective representation $U$ of $G$ with cocycle  
$\sigma\circ(q\times q)$, and  tensoring with an irreducible   
$\overline{\sigma}$-representation $W$ of $G/N$ gives an irreducible  
representation $U\otimes (W\circ q)$ of $G$ whose restriction to $N$ is  
a multiple $\pi\otimes 1$ of $\pi$.

\subsection*{Applications to compact Lie groups}\label{appl-cpt-lie}
When $\Gamma$ is a compact connected Lie group, Moore computed  
$H^2(\Gamma,\T)$ as follows. Let  $\widetilde \Gamma$ be the simply  
connected covering group of $\Gamma$; then the fundamental group  
$\pi_1(\Gamma)$ is isomorphic to a central subgroup of  
$\widetilde\Gamma$ and $\Gamma\cong\widetilde\Gamma/\pi_1(\Gamma)$. An  
inflation and restriction sequence identifies $H^2(\Gamma,\T)$ with the  
quotient of the dual group  
$\pi_1(\Gamma)^\wedge=\Hom(\pi_1(\Gamma),\T)$ by the image of the  
restriction map $\res:(\widetilde  
\Gamma)^\wedge\to\pi_1(\Gamma)^\wedge$ \cite[page 55]{Moore1}.

When $\Gamma=\T^n$, we have $\pi_1(\Gamma)=\Z^n$ and  
$\widetilde{\Gamma}=\R^n$, and the restriction map  
$\R^n=(\R^n)^\wedge\mapsto \T^n=(\Z^n)^\wedge$ is onto by duality. Thus  
$H^2(\T^n,\T)=0$. Theorem~\ref{thm-main} thus
implies that if $G/N\cong\T^n$, then every $G$-invariant irreducible unitary  
representation of $N$ extends to $G$.   Because representations of  
compact groups are direct sums of irreducible representations, this  
observation includes \cite[Corollary~3.5]{CKS}, and hence also  
\cite[Theorem~1.1]{CKS}.

For non-compact groups $G$, it is not clear how one might prove  
Corollary~\ref{cor-main} by reducing to the irreducible case: the  
analogue of the direct-sum decomposition would be a direct-integral  
decomposition, but not every unitary representation is a direct  
integral of irreducible representations. (See  \cite[\S10]{Mackey2} for  
a discussion of decomposing representations as direct integrals.)

\subsection*{The nonabelian duality approach}

If $\alpha:G\to \Aut(A)$ is a strongly continuous action of a locally  
compact group $G$ on a $C^*$-algebra $A$, a \emph{covariant
representation} of
$(A,G,\alpha)$ consists of a representation $\mu$ of $A$ and a unitary  
representation $U$ of $G$ such that
\[
\mu(\alpha_t(a))=U_t\mu(a)U_t^*\quad\text{for $a\in A$ and $t\in G$;}
\]
covariant representations can take values either in abstract  
$C^*$-algebras or in the concrete $C^*$-algebra $B(\H)$. The
\emph{crossed product}
$A\times_\alpha G$ is the $C^*$-algebra generated by a universal  
covariant representation in the multiplier algebra $M(A\times_\alpha  
G$) (see \cite{rae88} for details of what this means). The covariant
representations $(\mu,U)$ of $(A,G,\alpha)$ therefore give  
representations $\mu\times U$ of $A\times_\alpha G$, and
all representations of $A\times_\alpha G$ have this form. We shall be  
particularly interested in the actions $\lt:G\to
\Aut (C_0(G/N))$ and $\rt:G/N\to \Aut (C_0(G/N))$ defined by
\[
\lt_{s}(f)(uN)=f(s^{-1}uN)\quad \text{and}\quad
\rt_{tN}(f)(uN)=f(utN).
\]
The automorphisms $\rt_{tN}$ commute with the automorphisms $\lt_s$,  
and hence induce an action $\beta$ of
$G/N$ on the crossed product $C_0(G/N)\times_{\lt} G$.

If $\pi$ is a unitary representation of $N$, then the induced  
representation $\Ind \pi$ of $G$ acts in the completion
$\H(\Ind\pi)$ of
\[
\{\xi\in C_b(G,\H):\xi(tn)=\pi(n)^{-1}(\xi(t))\text{\ and\ }  
(tN\mapsto\|\xi(t)\|)\in C_c(G/N)\}\]
with respect to the inner product $(\xi\, |\, \eta)=\int_{G/N}  
(\xi(t)\, |\, \eta(t)) \, d(tN)$, according to the formula
$(\Ind\pi)_t(\xi)(r)=\xi(t^{-1}r)$.  (See, for example,  
\cite[page~296]{tfb}; because $N$ is normal there is a $G$-invariant  
measure on $G/N$, and we can take
the rho-function in the usual formula to be 1.)

Let
$M$ be the representation of
$C_0(G/N)$ by multiplication operators on
$\H(\Ind\pi)$, and note that
$(M,\Ind\pi)$ is a covariant representation of $(C_0(G/N), G,\lt)$.
The nonabelian duality approach to the extension problem yields the  
following theorem.

\begin{thm}\label{connect}
Suppose that $N$ is a closed normal subgroup of an amenable and  
second-countable locally compact group $G$, and suppose that $\pi$ is a  
  unitary
representation $\pi:N\to U(\H)$. Then $\pi$ extends to a  unitary  
representation of
$G$ if and only if there exists a unitary representation $Q:G/N\to  
U(\H(\Ind\pi))$ such that $(M\times\Ind\pi,Q)$ is a covariant  
representation
of $(C_0(G/N)\times_{\lt} G, G/N,\beta)$.
\end{thm}

\begin{proof}
The induction-restriction theory of \cite{BE}  says that $\pi$ is the  
restriction of a representation of $G$
if and only if $M\times \Ind\pi$ is induced, in a dual sense, from a  
representation of the group $C^*$-algebra $C^*(G)=\C\times G$.
To deduce this from \cite[Theorem~5.16]{BE}, we need to recall some  
ideas of nonabelian duality. The group $C^*$-algebra
$C^*(G)$ is generated by a universal unitary representation $\iota:G\to  
UM(C^*(G))$.
The comultiplication
$\delta:C^*(G)\to M(C^*(G)\otimes C^*(G))$ is the representation  
corresponding to the unitary representation
$\iota\otimes\iota$; it has a restriction $\delta|$ which is a coaction  
of $G/N$ on $C^*(G)$. Since $G$ is amenable, $C^*(G)$ coincides with  
the reduced group $C^*$-algebra $C^*_r(G)$, and hence we can apply  
results from \cite{BE} and \cite{Mansfield} concerning reduced crossed  
products. In particular, we can induce
representations from $C^*(G)$ to the coaction crossed product  
$C^*(G)\times_{\delta|} G/N$ by tensoring with a
$(C^*(G)\times_{\delta|}G/N)$--$C^*(G)$ bimodule $Y$ constructed by  
Mansfield \cite{Mansfield}; the resulting map on representations
is denoted by $Y\dashind$.

We   recall from \cite[Theorem~C.23]{tfb} that there is a Morita  
equivalence between $C_0(G/N)\times_{\lt} G$
and $C^*(N)$ which is implemented by an imprimitivity bimodule $X$; we  
denote by $X\dashind$ the corresponding map on representations. The  
algebras $C_0(G/N)\times_{\lt} G$ and
$C^*(G)\times_{\delta|}G/N$ have exactly the same covariant  
representations, and hence are isomorphic (see, for example,
\cite[Theorem~A.64]{BE}).
Thus we can view $X$ as a  $(C^*(G)\times_{\delta|} G/N)$--$C^*(N)$  
bimodule.
Theorem 5.16 of \cite{BE} (with $A=\C$ and $M=G$) says that, provided  
$G$ is amenable, we have a commutative diagram
\begin{diagram}
\Rep C^*(G) &\rTo(2.3,2.3)^{{\res}}   \Rep C^*(N) \\
        &\SE(1.9,1.9)_{{Y\dashind}} \hskip2.5cm  
\SW(1.9,1.9)_{{X\dashind}}\hskip1.5cm&\\
        &\Rep(C^*(G)\times_{\delta|}G/N)\hskip1cm&
\end{diagram}
  Since $X\dashind$ is a bijection, it follows that a representation  
$\pi$ of $C^*(N)$ extends to a representation of $C^*(G)$
if and only if $X\dashind\pi$ is in the range of $Y\dashind$.

To deduce Theorem~\ref{connect} from this, we have to make two  
observations. First, the representation
$X\dashind\pi$ of $C_0(G/N)\times_{\lt} G$ is equivalent to  
$M\times\Ind\pi$. To see this, note that the intertwining
unitary isomorphism
$W$ of $(X\otimes_{C^*(N)}\H, X\dashind\pi)$ onto  
$(\H(\Ind\pi),\Ind\pi)$ constructed in the proof of
\cite[Theorem~C.33]{tfb} carries the left action of $C_0(G/N)$ into  
$M$. Second, we recall from Mansfield's imprimitivity
theorem \cite[Theorem~28]{Mansfield} that a representation $\mu$ of  
$C^*(G)\times_{\delta|} G/N$ has the form $Y\dashind\rho$ if and only if
there is a unitary representation $Q$ of $G/N$ on $\H(\mu)$ such that  
$(\mu,Q)$ is covariant for the dual action
$(\delta|)^\wedge$ of $G/N$. Since \cite[Theorem~A.64]{BE} also says  
that the isomorphism of $C_0(G/N)\times_{\lt}G$ onto
$C^*(G)\times_{\delta|} G/N$ carries the action $\beta$ into  
$(\delta|)^\wedge$, the result follows.
\end{proof}

\begin{remark}\label{remark}
Comparing Theorem~\ref{connect} with Theorem~\ref{thm-main}, it is  
natural to ask what happened to the hypothesis ``$\pi$ is  
$G$-invariant''.
Suppose $\pi$ is $G$-invariant, so that there exist unitary operators  
$W_s$ on $\H$ such that $W_s\pi(n)W_s^*=\pi(sns^{-1})$. Then
\begin{equation*}\label{eq-U}
U_s(\xi)(t)=W_s(\xi(ts))\Delta(sN)^{1/2}
\end{equation*}
defines a unitary operator $U_s$ on $\H(\Ind\pi)$ which intertwines the  
covariant representations $(M,\Ind\pi)$ and $(M\circ \rt_{sN},  
\Ind\pi)$.  So $R_{sN}:=U_{c(sN)}$ defines a map $R:G/N\to  
U(\H(\Ind\pi))$ which formally satisfies the covariance relations but  
is not necessarily a representation.

Our original extension problem for a $G$-invariant representation  
$\pi:N\to U(\H)$ therefore reduces to:
\begin{itemize}
\item[] Given a representation $\phi$ of $C^*(G)\times_{\delta|}G/N$  
such that $\phi\circ(\delta|)^{\wedge}_{sN}$ is equivalent to $\phi$  
for every $sN\in G/N$, is there a representation $Q$ of $G/N$ such that  
$(\phi,Q)$ is covariant for $(C^*(G)\times_{\delta|}G/N, G/N,  
(\delta|)^\wedge)$?
\end{itemize}
Since there are by hypothesis unitary operators $R_{sN}$ such that  
$\phi\circ(\delta|)^{\wedge}_{sN}=\Ad R_{sN}\circ\phi$, we can repeat  
the analysis  of Theorem~\ref{thm-main} to see that there is a twisted  
action $(\beta,\omega)$ of $G/N$ on the commutant of the range of  
$\phi$, such that $(\beta,\omega)$ is exterior equivalent to an  
ordinary action if and only if we can adjust the $R_{sN}$ to obtain the  
required representation $Q$.
Thus $\beta_{sN}=\Ad R_{sN}$ and, for $\xi\in\H(\Ind\pi)$,
\begin{align*}
\omega(rN,sN)(\xi)(t)&=R_{rN}R_{sN}R_{rsN}^*(\xi)(t)=
U_{c(rN)}U_{c(sN)}U_{c(rsN)}^*(\xi)(t)\\
&=W_{c(rN)}W_{c(sN)}W_{c(rsN)}^*(\xi(tc(rN)c(sN)c(rsN)^{-1}))\\
&=W_{c(rN)}W_{c(sN)}W_{c(rsN)}^*\pi(c(rN)c(sN)c(rsN)^{-1})^{- 
1}(\xi(t))\\
&=W_{c(rN)}W_{c(sN)}W_{c(rsN)}^*\pi(c(rsN)c(sN)^{-1}c(rN)^{-1})(\xi(t)).
\end{align*}

We claim that the obstruction $(\beta,\omega)$ is essentially the same  
as the obstruction $(\alpha,\sigma)$ to extending $\pi$ from  
Theorem~\ref{thm-main}.  To see this, we first  identify  
$\pi(N)'$ with $\phi(C^*(G)\times_{\delta|}G/N)'$ when  
$\phi=M\times\Ind\pi$. If $T\in\pi(N)'$, then the formula $1\otimes  
T(\xi)(t)=T(\xi(t))$ defines an operator in  
$\phi(C^*(G)\times_{\delta|}G/N)'$.  When we view $\H(\Ind\pi)$ as  
$X\otimes_{C^*(N)}\H$, then we recover $\H$ as $\widetilde  
X\otimes_{C^*(G)\times G/N}(X\otimes_{C^*(N)}\H)$, where $\widetilde X$  
is the dual imprimitivity bimodule, and the natural isomorphism of $\H$  
onto $\widetilde X\otimes_{C^*(G)\times G/N}(X\otimes_{C^*(N)}\H)$  
takes $T$ to $1\otimes 1\otimes T$.  Thus $T\mapsto 1\otimes T$ is an  
isomorphism of $\pi(N)'$ onto $\phi(C^*(G)\times_{\delta|}G/N)'$.

With $V$ as in Equation~\eqref{eq-V}, the cocycle $\sigma$ in the  
twisted action $(\alpha,\sigma)$ satisfies
\begin{align*}
\sigma(rN,sN)&=V_rV_sV_{rs}^*\\
&=W_{c(rN)}\pi(c(rN)^{-1}r)W_{c(sN)}\pi(c(sN)^{-1}s)
\pi(c(rsN)^{-1}rs)^{-1}W_{c(rsN)}^*\\
&=W_{c(rN)}W_{c(sN)}\pi(c(sN)^{-1}c(rN)^{-1}rc(sN))\pi(c(sN)^{-1}s)
\pi(s^{-1}r^{-1}c(rsN))W_{c(rsN)}^*\\
&=W_{c(rN)}W_{c(sN)}\pi(c(sN)^{-1}c(rN)^{-1}c(rsN))W_{c(rsN)}^*\\
&=W_{c(rN)}W_{c(sN)}W_{c(rsN)}^*\pi(c(rsN)c(sN)^{-1}c(rN)^{-1}).
\end{align*}
Thus with this choice of $R_{sN}$, we have
$\omega(rN,sN)=1\otimes\sigma(rN,sN)$, and for $T\in\pi(N)'$,
\begin{align*}
\beta_{sN}(1\otimes T)(\xi)(t)&=R_{sN}(1\otimes T)R_{sN}^*(\xi)(t)\\
&=W_{c(sN)}TW_{c(sN)}^*(\xi(t))\\
&=(1\otimes V_sTV_s^*)(\xi)(t)\\
&=(1\otimes\alpha_s(T))(\xi)(t).
\end{align*}
So the isomorphism of $\pi(N)'$ onto $\phi(C^*(G)\times_{\delta|}G/N)'$  
carries $(\alpha, \sigma)$ into the twisted action $(\beta,\omega)$  
which obstructs the existence  of $Q$.
Thus, reassuringly, the cohomological obstruction to finding $Q$ is  
identical to the obstruction to extending $\pi$. 
\end{remark}


\begin{thebibliography}{99}
\bibitem{CKS} J.W. Cho,  M.K. Kim and  D.Y Suh, \emph{On extensions of  
representations for compact Lie groups}, J. Pure Appl. Algebra  
\textbf{178} (2003), 245--254.

\bibitem{Clifford}  A.H. Clifford, \emph{Representations induced in an  
invariant subgroup}, Ann. of Math. \textbf{38} (1937), 533--550.

\bibitem{BE} S. Echterhoff, S. Kaliszewski, J. Quigg and I.
   Raeburn, \emph{A categorical approach to imprimitivity theorems for  
$C^*$-dynamical systems}, preprint, arXiv:math.OA/0205322.

\bibitem{Haagerup}  U. Haagerup, \emph{The standard form of von Neumann  
algebras}, Math Scand. \textbf{37} (1975), 271--283.

\bibitem{KR} R.V. Kadison  and  J.R. Ringrose, Fundamentals of the  
theory of operator algebras \textrm{II}, Graduate Studies in  
Mathematics, vol. 16,
Amer. Math. Soc., Providence,  1997.

\bibitem{Mackey2}  G.W. Mackey, \emph{Borel structure in groups and  
their duals}, Trans. Amer. Math. Soc. \textbf{85} (1957), 134--165.

\bibitem{Mackey}  G.W. Mackey, \emph{Unitary representations of group  
extensions.
\textrm{I}}, Acta Math \textbf{99} (1958), 265--311.

\bibitem{Mansfield}
  K. Mansfield, \emph{Induced representations of crossed products by  
coactions}, J. Funct. Anal. \textbf{97} (1991), 112--161.

\bibitem{Moore1}  C.C. Moore, \emph{Extensions and low dimensional  
cohomology theory of locally compact groups. \textrm{I}}, Trans. Amer.  
Math. Soc. \textbf{113} (1964), 40--63.

\bibitem{Moore2}  C.C. Moore, \emph{Extensions and low dimensional  
cohomology theory of locally compact groups. \textrm{II}}, Trans. Amer.  
Math. Soc. \textbf{113} (1964), 64--86.

\bibitem{Moore3}  C.C. Moore, \emph{Group extensions and cohomology for  
locally compact groups. \textrm{III}}, Trans. Amer. Math. Soc.  
\textbf{221} (1976), 1--33.

\bibitem{Moore4}  C.C. Moore, \emph{Group extensions and cohomology for  
locally compact groups. \textrm{IV}}, Trans. Amer. Math. Soc.  
\textbf{221} (1976), 35--58.

\bibitem{PR} J.A. Packer and I. Raeburn, \emph{Twisted crossed products  
of $C^*$-algebras}, Math. Proc. Camb. Phil. Soc. \textbf{106} (1989),  
293--311.

\bibitem{rae88} I. Raeburn, \emph{On crossed products and Takai  
duality},
  Proc. Edinburgh Math. Soc. \textbf{31} (1988), 321--330.

\bibitem{tfb}
I. Raeburn and D.P. Williams, Morita Equivalence and
Continuous-Trace $C^*$-Algebras, Math. Surveys and Monographs, vol.
60, Amer. Math. Soc., Providence, 1998.

\bibitem{Colin1}  C.E. Sutherland, \emph{Cohomology and extensions of  
von Neumann algebras. \textrm{I}},  Publ. Res. Inst. Math. Sci.  
\textbf{16} (1980), 105--133.


\bibitem{Colin2}
  C.E. Sutherland, \emph{Cohomology and extensions of von Neumann  
algebras. \textrm{II}},
Publ. Res. Inst. Math. Sci. \textbf{16} (1980), 135--174.

\end{thebibliography}
\end{document}